\documentclass{amsart}        

\usepackage{amsthm}   

\newtheorem{cor}{Corollary}
\newtheorem{cnj}[cor]{Conjecture}

\newcommand{\dfn}{\textbf} 

\newcommand{\ie}                {\emph{i.e.,}}

\DeclareMathOperator{\susp}{susp}

\newdimen\bboxht  \newdimen\bboxwd  
\def\bbox#1#2#3#4#5{%
\bboxwd=-#1 \advance\bboxwd by#3    
\bboxht=-#2 \advance\bboxht by#4    
\setbox0=\hbox{\hskip-#1\lower#2\hbox{#5}}%
\ht0=\bboxht \wd0=\bboxwd \dp0=0pt \box0 }

\newcommand{\PSh}[1]{\special{! #1}}
\newcommand{\PS}[1]{\special{" #1}}
\newcommand{\PSbox}[5]{\bbox{#1}{#2}{#3}{#4}{\PS{#5}}}

\PSh{
/r 2.1 def      /lw 1.3 def     /offset 12 def    /fs 9.5 def
2 setlinejoin
/Circ { newpath r 0 360 arc eofill } def
/Vertex { 2 copy Circ newpath moveto
dup  cos offset mul
exch sin offset mul
/Times-Roman findfont fs scalefont setfont
rmoveto -3 -3 rmoveto 
show } def
/IVertex { 2 copy Circ newpath moveto
dup  cos offset mul
exch sin offset mul
/Times-Italic findfont fs scalefont setfont
rmoveto -3 -3 rmoveto 
show } def
/Path { lw setlinewidth 
3 1 roll newpath moveto -1 add { lineto } repeat stroke } def
/Poly { lw setlinewidth
3 1 roll newpath moveto -1 add { lineto } repeat closepath
gsave setgray eofill grestore stroke } def
/base 40 def    
/h { base 3 sqrt mul } def
/v1 { 0 0 } def
/v2 { base h } def
/v3 { base 2 mul 0 } def
/v4 { base 3 mul h } def
/v5 { base h 1.5 mul } def
/v6 { base h 2   mul } def
/shade .85 def  /White 1 def
}

\begin{document}

\title{David Gale's Subset Take-away Game}

\author{J. Daniel Christensen}
\address{Department of Mathematics\\ M.I.T.\\ Cambridge, MA 02139} 
\email{jdchrist@math.mit.edu}

\author{Mark Tilford}
\address{MSC\#936 Caltech\\Pasadena, CA 91126}
\email{tilford@cco.caltech.edu}

\thanks{This version is more up to date than the published version.  
See the last page for details.}

\maketitle

  \dfn{Subset take-away} is a two-player game involving a fixed
finite set $A$.
Players alternately choose proper, non-empty subsets of $A$, 
with the condition that one may not name a set containing a set
that was named earlier.  A player who is unable to move loses.  
For example, if $A = \{ 1 \}$, then there are no
legal moves and the second player wins.
If $A = \{1, 2\}$, then the only legal moves are $\{ 1 \}$ and $\{ 2 \}$.
Each is a good reply to the other, and so once again the second
player wins.
The first interesting case is when $A = \{1, 2, 3\}$.  In response
to any first move, the second player may choose the complementary set.
This produces a position equivalent to the starting position when
$A = \{1, 2\}$ and thus leads to a win for the second player.
With increasing patience, the reader may enjoy verifying that
when $A$ has fewer than $6$ elements, the game is a second
player win.
Indeed, David Gale, whom we understand deserves credit for this game, 
made the following conjecture \cite{gu:upcg}.

\begin{cnj}
  Subset take-away is always a second player win.
\end{cnj}

It was pointed out to us by
Richard Stanley that the collection of legal moves at any given
time forms an abstract simplicial complex: if $X$ is a subset of $A$ that
is legal, then every non-empty subset of $X$ is also legal.
To translate this into geometry, view a $k$-subset 
as a $(k-1)$-simplex. 
If $A = \{1, 2, \cdots, n\}$, then the starting position
corresponds to the boundary of the $(n-1)$-simplex.
A move in this formulation of the game consists of
choosing a simplex of any dimension, and erasing its interior
as well as all higher-dimensional simplices having it as a face.
For example, the starting position for the $n = 4$ game is a hollow
tetrahedron, and after the move $\{1, 4\}$ (which we write 14 for brevity)
we have the position
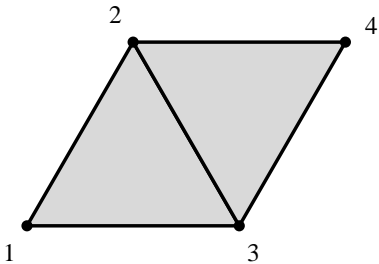
\[
\PSbox{-23bp}{-23bp}{130bp}{130bp}{
shade v1 v2 v3 3 Poly
shade v2 v3 v4 3 Poly
(1) -120 v1 Vertex
(2) 120 v2 Vertex
(3) -60 v3 Vertex
(4)  30 v4 Vertex
}
\]
in which the legal moves are 123 and 234 along with their
non-empty subsets.
Now one could remove the 2-simplex 234 leaving
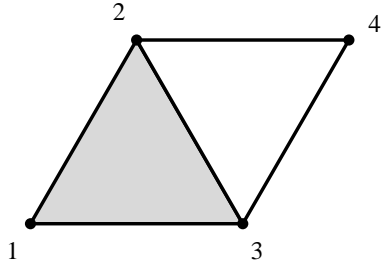
\[
\PSbox{-23bp}{-23bp}{130bp}{100bp}{
shade v1 v2 v3 3 Poly
White v2 v3 v4 3 Poly
(1) -120 v1 Vertex
(2) 120 v2 Vertex
(3) -60 v3 Vertex
(4)  30 v4 Vertex
} . 
\]
And if the next player removed the vertex 3, the new position would
be
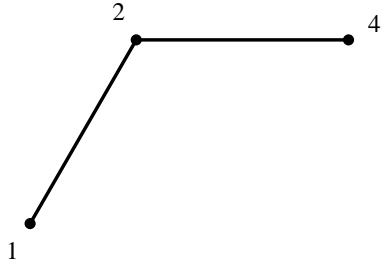
\[
\PSbox{-23bp}{-23bp}{130bp}{100bp}{
v1 v2 v4 3 Path
(1) -120 v1 Vertex
(2) 120 v2 Vertex
(4)  30 v4 Vertex
} \quad ,
\]
with five legal moves.  A winning move here is to erase the vertex 2,
leaving a game equivalent to the starting position for $n = 2$.
The geometric formulation is quite helpful, and we will use this
language frequently.  

We next describe a technique that allows one to reduce the size of
a position without changing its win/loss value.
Using this method one can show in just a
couple of minutes that the $n = 5$ game is a second player win.  Also,
the $n = 6$ game can be reduced to a simpler game that we were able to
analyze with a computer.  Assuming the correctness of our program,
we can assert that the $n = 6$ game is also a second player win.

\section{The reduction technique}

  We should first mention that while the game can be thought of
geometrically, it is still a combinatorial game in the sense that the
choice of triangulation matters.  For example, the simplicial
complexes
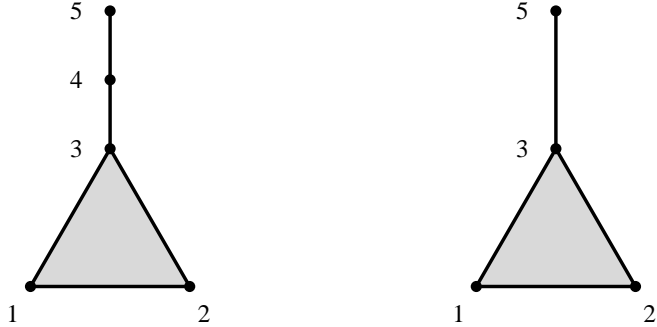
\[
\PSbox{-18bp}{-18bp}{78bp}{119bp}{
/base 30 def
shade v1 v2 v3 3 Poly
v2 v5 v6 3 Path
(1) -120 v1 Vertex
(3) 180 v2 Vertex
(2)  -60 v3 Vertex
(4) 180 v5 Vertex
(5) 180 v6 Vertex
}
\hspace{1.0in}
\PSbox{-18bp}{-18bp}{78bp}{119bp}{
/base 30 def
shade v1 v2 v3 3 Poly
v2 v6 2 Path
(1) -120 v1 Vertex
(3) 180 v2 Vertex
(2)  -60 v3 Vertex
(5) 180 v6 Vertex
}
\]
are different triangulations of the same space, 
but the first is a first player win, while
the second is a second player win.  
A winning move in the first game is to erase the vertex 3, leaving
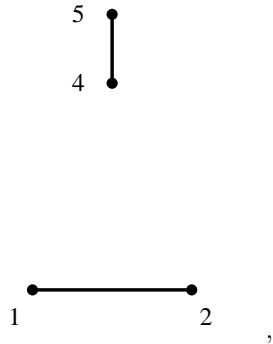
\[
\PSbox{-18bp}{-18bp}{78bp}{119bp}{
/base 30 def
v1 v3 2 Path
v5 v6 2 Path
(1) -120 v1 Vertex
(2)  -60 v3 Vertex
(4) 180 v5 Vertex
(5) 180 v6 Vertex
}
\quad , \]
which is a second player win by symmetry.
In the second game, the second player responds to the first
according to the following strategy:
\newcommand{\lra}{\ \longleftrightarrow\ }
\newcommand{\qq}{\qquad}
\[ 5 \lra 123 \qq 3 \lra \text{$1$ or $2$} \qq 12 \lra 35 \qq 13 \lra 23 . \]

  Now, despite the fact that the game seems to have little to do with
topology, Keith Orpen suggested that we think about suspension
since the conjecture deals with spheres, and the suspension of a
$k$-sphere is a $(k+1)$-sphere.  If $X$ is a collection of subsets of
$A$ which form a simplicial complex, then the \dfn{suspension} $\susp X$
of $X$ is the simplicial complex consisting of the following 
subsets of $A \cup \{ x, y \}$, where $x$ and $y$ are new
vertices:  $\susp X$ contains the sets in $X$, the sets $a \cup \{x\}$
and $a \cup \{y\}$ for each set $a$ in $X$, and the sets $\{x\}$ and $\{y\}$.
For example, the suspension of the interval is a diamond:
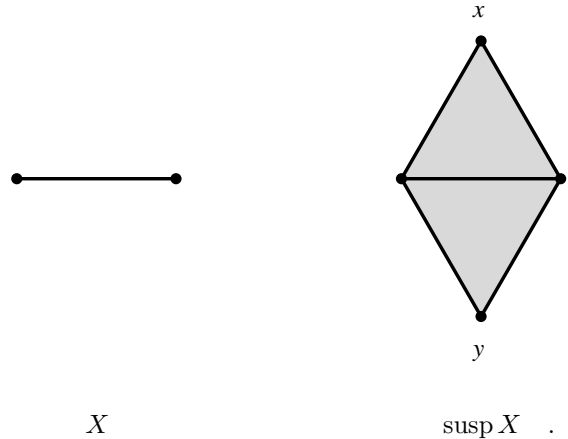
\PSh{
/base2 30 def 
/a { 0 0 } def
/b { base2 2 mul 0 } def
/x { base2 base2 3 sqrt mul } def
/y { x neg } def
/c { base2 4 mul 0 } def
/d { base2 6 mul 0 } def
}
\[ 
\begin{array}{c}
  \PSbox{-1bp}{-84bp}{62bp}{69bp}{
  a b 2 Path  a Circ  b Circ} \\ X
\end{array}
\hspace{1.0in}
\begin{array}{c}
  \PSbox{-1bp}{-84bp}{62bp}{69bp}{
  shade a b x 3 Poly
  shade a b y 3 Poly
  a Circ   b Circ
  (x) 90   x IVertex
  (y) -90  y IVertex
} \\
  \susp X \rlap{\quad .}
\end{array}
\]
It turns out that $X$ and $\susp X$ always have the same win/loss value,
and so a position of the form $\susp X$ can be \emph{reduced}.
However, if $X$ is the $k$-sphere triangulated as the boundary of
the $(k+1)$-simplex, then $\susp X$ is topologically a $(k+1)$-sphere
but is \emph{not} triangulated as the boundary of the $(k+2)$-simplex.
Therefore this reduction can't be used directly to attack the conjecture.

There is a more general type of reduction that is possible.  
We say that a pair $(x,y)$ of vertices in a simplicial complex $X$
is a \dfn{binary star} if
\begin{itemize}
  \item there is no edge connecting $x$ and $y$ (\ie\ $\{x, y\}$ is not present), 
  \item for each set $a$ in $X$, if $a$ contains $x$ then there 
        is a set $b$ which is the same as $a$ except that it 
        contains $y$ in place of $x$, and
  \item the same as the previous item with $x$ and $y$ interchanged.
\end{itemize}
When this is the case, the vertices $x$ and $y$ can be removed (along with
the simplices containing them) without changing who wins.  
Indeed, if player A
can win in the reduced game, then this player can win in the larger
game by using the following strategy:
when B makes a move involving $x$ or $y$, A makes the corresponding
move with $x$ and $y$ interchanged;  when B makes a move not involving
$x$ or $y$, A replies with the winning response in the reduced game.


  Here's an example of the power of binary star reduction.  Consider
the starting position when $A = \{1, 2, \cdots, n\}$.
If an edge is removed, then
its endpoints form a binary star in the resulting position.
Performing reduction produces a (filled in) $(n-3)$-simplex.  It is easy
to see that this position is a first player
win, since the first player can ``pass'' by removing the interior if
no other move wins, thereby leaving the second player with a losing
position.  Thus we have proved that choosing an edge from the starting
position is a losing move.

  Similarly, one can show that moves of size $1$, $n-2$ and $n-1$ are all
bad in the standard starting position of size $n$, assuming Gale's
conjecture for smaller starting positions.  
(By a move of size $1$ we mean a vertex, etc.)
When $n = 5$ these are the only moves available, so Gale's conjecture
holds in this case.  
A reader who verified this without binary star reduction
will realize how much effort we have saved.
When $n = 6$ we only have to analyze the situation after 
playing a move of size $3$.
Based on our experience with smaller games we guessed
that the complementary set is a winning response, and the computer
verified this.  In this way we verified Gale's conjecture for $n = 6$
and were led to the following conjecture.
\begin{cnj}
  A winning response to an opening move is to 
  play the complementary move.
\end{cnj}

\section{A counterexample}

\PSh{ 
/base3 25 def 
/c { 0 0 } def
/d { base3 0 } def
/e { base3 2 mul 0 } def
/f { base3 3 mul 0 } def
}
The game
\[ \PSbox{-2bp}{-2bp}{2bp}{2bp}{c Circ} \]
is clearly a first player win.  Also, each of the games homeomorphic
to the 1-simplex, namely
\[ 
\PSbox{-2bp}{-2bp}{29bp}{2bp}{c d 2 Path c Circ d Circ}, \quad
\PSbox{-2bp}{-2bp}{54bp}{2bp}{c d e 3 Path c Circ d Circ e Circ}, \quad
\PSbox{-2bp}{-2bp}{79bp}{2bp}{c d e f 4 Path c Circ d Circ e Circ f Circ},\quad
\ldots
\]
is a first player win since the first player can choose the central
simplex and then play by symmetry.  
One might hope that any triangulation of the $k$-simplex is
a first player win.
(Indeed, this implies Gale's conjecture.) 
However, the following is a counterexample:
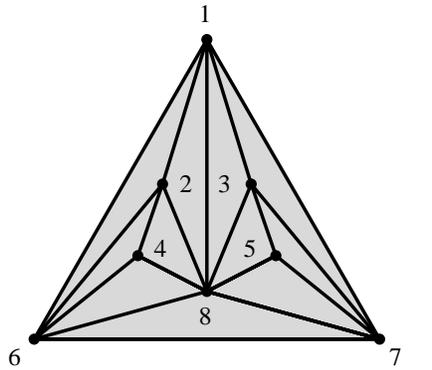
\[ 
\PSbox{-12bp}{-12bp}{140bp}{129bp}{
  /offset 9.5 def
  2 setlinejoin
  65 0 translate
  /w1 {  0   113 } def
  /w2 { -16.7 58.5 } def
  /w3 {  16.7 58.5 } def
  /w4 { -26   31.5 } def
  /w5 {  26   31.5 } def
  /w6 { -65.1 0   } def
  /w7 {  65.1 0   } def
  /w8 {  0   18  } def
  shade w1 w6 w7 3 Poly
  w1 w2 w6 w4 w2 w8 w4 w8 w6 w7 w8 w3 w5 w8 w1 w3 w7 w5 w8 w7 20 Path
  (1)  90 w1 Vertex
  (2)   0 w2 Vertex
  (3) 180 w3 Vertex
  (4)  15 w4 Vertex
  (5) 165 w5 Vertex
  (6) 225 w6 Vertex
  (7) -45 w7 Vertex
  (8) -90 w8 Vertex
} \quad . \]
This triangulation of the 2-simplex is a second player win.
Notice that by adding the 2-simplex 167 to the above picture,
one gets a triangulation of the 2-sphere which is a first player win.
We wonder whether there are simpler triangulations of the
2-simplex and the 2-sphere that have these win/loss values.


\newpage

\noindent
\textbf{\large This version has been updated compared to the published version:}

\medskip

\noindent
The published version appeared in
\emph{American Mathematical Monthly} \textbf{104} (1997), 762--766.
Since then, the following changes have been made:

\begin{itemize}
\item In the proof that choosing an edge is a losing move, ``$(n-2)$-simplex''
      has been corrected to ``$(n-3)$-simplex''.
\item The reference to \cite{gu:upcg} has been updated.
\item Minor formatting changes have been made.
\end{itemize}

The contact information is out of date.
For questions and comments, contact Dan Christensen at \texttt{jdc@uwo.ca}.


\begin{thebibliography}{1}

\bibitem{gu:upcg}
Richard~K. Guy.
\newblock Unsolved problems in combinatorial games.
\newblock \emph{Games of No Chance}, (R.\ J.\ Nowakowski ed.)
\newblock MSRI Publications \textbf{29}, Cambridge University Press, 1996,
          pp.\ 475--491.

\end{thebibliography}
\end{document}